\documentclass[11pt]{article}
\usepackage{amscd}
\usepackage{amssymb,amsfonts,amsmath,psfrag,amscd,cite}
\usepackage{amsmath}
  \usepackage{paralist}
  \usepackage{graphics}
  \usepackage{epsfig}
  \usepackage[colorlinks=true]{hyperref}
  \usepackage[all]{xy}
\newtheorem{theo}{Theorem}

\newtheorem{defi}[theo]{Definition}

\newtheorem{rema}[theo]{Remark}

\makeatletter \@addtoreset{equation}{section}
\@addtoreset{theo}{section} \makeatother

\begin{document}


\date{}
\title{The Geometrical Structure of Phase Space of
\\ the Controlled Hamiltonian System with Symmetry}
\author{Hong Wang \\
School of Mathematical Sciences and LPMC,\\
Nankai University, Tianjin 300071, P.R.China\\
E-mail: hongwang@math.nankai.edu.cn}
\date{\emph{ Dedicated to My Advisor---Professor Hesheng Hu\\
   on the Occasion of Her 90th Birthday}\\
   November 1, 2020 } \maketitle

{\bf Abstract:} In this paper, from the viewpoint of completeness of
Marsden-Weinstein reduction, we illustrate how to give the
definitions of a controlled Hamiltonian (CH) system and a reducible
controlled Hamiltonian system with symmetry; and how to describe
the dynamics of a CH system and the controlled Hamiltonian
equivalence; as well as how to give the regular point reduction and
the regular orbit reduction for a CH system with
symmetry, by analyzing carefully the geometrical and topological
structures of the phase space and the reduced phase space of the
corresponding Hamiltonian system. We also introduce
briefly some recent developments in the study of reduction
theory for the CH systems with symmetries and applications.
These research work reveal the deeply internal
relationships of the geometrical structures of phase spaces, the dynamical
vector fields and controls of the CH systems.\\




\tableofcontents

\section{Introduction}

The reduction theory for a mechanical system with symmetry is an
important subject and it is widely studied in the theory of
mathematics and mechanics, as well as applications. The main goal of
reduction theory in mechanics is to use conservation laws and the
associated symmetries to reduce the number of dimensions of a
mechanical system required to be described. So, such reduction
theory is regarded as a useful tool for simplifying and studying
concrete mechanical systems. Over forty years ago, the regular
symplectic reduction for the Hamiltonian system with symmetry and
coadjoint equivariant momentum map was set up by famous professors
Jerrold E. Marsden and Alan Weinstein, which is called
Marsden-Weinstein reduction, and great developments have been
obtained around the work in the theoretical study and applications
of mathematics, mechanics and physics; see
Abraham and Marsden \cite{abma78}, Abraham et al.
\cite{abmara88}, Arnold \cite{ar89}, Libermann and Marle
\cite{lima87}, Marsden \cite{ma92}, Marsden et al.
\cite{mamiorpera07, mamora90}, Marsden and Perlmutter \cite{mape00},
Marsden and Ratiu \cite{mara99}, Marsden and
Weinstein \cite{mawe74}, Meyer \cite{me73},
Nijmeijer and Van der Schaft \cite {nivds90} and
Ortega and Ratiu \cite{orra04}.\\

Recently, in Marsden et al.\cite{mawazh10}, the authors found that
the symplectic reduced space of a Hamiltonian system
defined on the cotangent bundle of a configuration manifold may not be a
cotangent bundle, and hence the set of Hamiltonian systems with
symmetries on the cotangent bundle is not complete under the
Marsden-Weinstein reduction. This is a serious problem. If we define
directly a controlled Hamiltonian system with symmetry on a
cotangent bundle, then it is possible that the Marsden-Weinstein
reduced system may not have definition. The study of completeness
of Hamiltonian reductions for a Hamiltonian system with
symmetry is related to the geometrical and topological structures of
Lie group, configuration manifold and its cotangent bundle, as well
as the action way of Lie group on the configuration manifold and its
cotangent bundle. In order to define the CH system and set up
the various perfect reduction
theory for the CH systems, we have to give the
precise analysis of geometrical and topological structures of the
phase spaces and the reduced phase spaces of various CH systems.\\

A controlled Hamiltonian (CH) system is a
Hamiltonian system with external force and control. In general, a CH
system under the actions of external force and control is not
Hamiltonian, however, it is a dynamical system closely related to a
Hamiltonian system, and it can be explored and studied by extending
the methods for external force and control in the study of Hamiltonian systems.
Thus, we can emphasize explicitly the impact of external force
and control in the study of CH systems. For example,
in order to describe the feedback control law to
modify the structures of CH system and the reduced CH system,
we introduce the notions of CH-equivalence, RpCH-equivalence and
RoCH-equivalence and so on.
In this paper, at first, from the viewpoint of completeness of
Marsden-Weinstein reduction, we shall illustrate how to give the
definitions of a CH system and a reducible CH system;
Next, we illustrate how to describe the dynamics of a
CH system and the controlled Hamiltonian
equivalence; The third, we give the regular point reduction and
the regular orbit reduction for a CH system with
symmetry and momentum map, by analyzing carefully the geometrical and topological
structures of the phase space and the reduced phase space of the
corresponding Hamiltonian system. Finally, we also introduce
briefly some recent developments in the study of reduction
theory for the CH systems with symmetries and applications.
These research works not only gave a
variety of reduction methods for the CH systems, but
also showed a variety of relationships of the
controlled Hamiltonian equivalences and the regular reducible
controlled Hamiltonian equivalences.
Moreover, we also state that the geometric constraint conditions of
canonical symplectic form and regular reduced symplectic forms for the
dynamical vector fields of an RCH system and its regular
reduced systems can be derive precisely,
which are called the Type I and Type II of Hamilton-Jacobi equations.\\

\section{Marsden-Weinstein Reduction on a Cotangent Bundle}

It is well-known that, in mechanics, the phase space of a
Hamiltonian system is very often the cotangent bundle $T^*Q$ of a
configuration manifold $Q$, and the reduction theory on the cotangent
bundle of a configuration manifold is a very important special case
of general symplectic reduction theory. In the following we first
give the Marsden-Weinstein reduction for a Hamiltonian system
with symmetry on the cotangent bundle of a smooth
configuration manifold with canonical symplectic structure,
see Abraham and Marsden \cite{abma78} and
Marsden and Weinstein \cite{mawe74}.\\

Let $Q$ be a smooth manifold and $TQ$ the tangent bundle, $T^* Q$
the cotangent bundle with a canonical symplectic form $\omega_0$.
Assume that $\Phi: G\times Q \rightarrow Q$ is a left smooth action
of a Lie group $G$ on the manifold $Q$. The cotangent lift is the
action of $G$ on $T^\ast Q$, $\Phi^{T^\ast}:G\times T^\ast
Q\rightarrow T^\ast Q$ given by $g\cdot
\alpha_q=(T\Phi_{g^{-1}})^\ast\cdot \alpha_q,\;\forall\;\alpha_q\in
T^\ast_qQ,\; q\in Q$. The cotangent lift of any proper (resp. free)
$G$-action is proper (resp. free). Assume that the cotangent lift action is
symplectic with respect to the canonical symplectic form $\omega_0$,
and has an $\operatorname{Ad}^\ast$-equivariant momentum map
$\mathbf{J}:T^\ast Q\to \mathfrak{g}^\ast$ given by
$<\mathbf{J}(\alpha_q),\xi>=\alpha_q(\xi_Q(q)), $ where $\xi\in
\mathfrak{g}$, $\xi_Q(q)$ is the value of the infinitesimal
generator $\xi_Q$ of the $G$-action at $q\in Q$, $<,>:
\mathfrak{g}^\ast \times \mathfrak{g}\rightarrow \mathbb{R}$ is the
duality pairing between the dual $\mathfrak{g}^\ast $ and
$\mathfrak{g}$. Assume that $\mu\in \mathfrak{g}^\ast $ is a regular
value of the momentum map $\mathbf{J}$, and $G_\mu=\{g\in
G|\operatorname{Ad}_g^\ast \mu=\mu \}$ is the isotropy subgroup of
the coadjoint $G$-action at the point $\mu$. From the
Marsden-Weinstein reduction, we know that the reduced space
$((T^*Q)_\mu, \omega_\mu)$ is a symplectic manifold.\\

In the following we give further a precise analysis of the
geometrical structure of the symplectic reduced space of $T^* Q$.
From Marsden and Perlmutter \cite{mape00} and Marsden et al.
\cite{mamiorpera07}, we know that the classification of symplectic
reduced space of the cotangent bundle $T^* Q$ as follows. (1) If
$\mu=0$, the symplectic reduced space of cotangent bundle $T^\ast Q$
at $\mu=0$ is given by $((T^\ast Q)_\mu, \omega_\mu)= (T^\ast(Q/G),
\omega_0)$, where $\omega_0$ is the canonical symplectic form of the
reduced cotangent bundle $T^\ast (Q/G)$. Thus, the symplectic
reduced space $((T^\ast Q)_\mu, \omega_\mu)$ at $\mu=0$ is a
symplectic vector bundle; (2) If $\mu\neq0$, and $G$ is Abelian,
then $G_\mu=G$, in this case the regular point symplectic reduced
space $((T^*Q)_\mu, \omega_\mu)$ is symplectically diffeomorphic to
symplectic vector bundle $(T^\ast (Q/G), \omega_0-B_\mu)$, where
$B_\mu$ is a magnetic term; (3) If $\mu\neq0$, and $G$ is not
Abelian and $G_\mu\neq G$, in this case the regular point symplectic
reduced space $((T^*Q)_\mu, \omega_\mu)$ is symplectically
diffeomorphic to a symplectic fiber bundle over $T^\ast (Q/G_\mu)$
with fiber to be the coadjoint orbit $\mathcal{O}_\mu$, see the
cotangent bundle reduction theorem---bundle version, also see
Marsden and Perlmutter \cite{mape00} and Marsden et al.
\cite{mamiorpera07}.\\

Thus, from the above discussion, we know that
the symplectic reduced space of a Hamiltonian system
defined on the cotangent bundle of a configuration manifold may not be a
cotangent bundle. Therefore, the symplectic reduced system of a
Hamiltonian system with symmetry defined on the cotangent bundle
$T^*Q$ may not be a Hamiltonian system on a cotangent bundle, that
is, the set of Hamiltonian systems with symmetries on the cotangent
bundle is not complete under the Marsden-Weinstein reduction.

\section{The Definition of Controlled Hamiltonian System}

In 2005, we hope to study the Hamiltonian system with control from
geometrical viewpoint. I and my students read the following two
papers of Professor Jerrold E. Marsden and his students on a seminar of Nankai University,
see Chang et al. \cite{chbllemawo02, chma04}. We found that
there are some serious wrong of rigor for the definitions
of controlled Lagrangian (CL) system, controlled Hamiltonian (CH)
system and its reduced CH systems, as well as
CH-equivalence and the reduced CH-equivalence
in this two papers.\\

There are the following four aspects of these wrong
for CL system, CH system and the reduced CH system:\\

(1) The authors define CL system by using a wrong expression.
In fact, the CL system is defined in \cite{chbllemawo02,chma04},
by using the following expression
\begin{equation}
\mathcal{EL}(L)(q,\dot{q},\ddot{q})=F(q,\dot{q})+u(q,\dot{q}),
\end{equation}
where $\mathcal{EL}$ is the Euler-Lagrange operator and the bundle maps
$\mathcal{EL}(L): T^{(2)}Q \rightarrow T^*Q, $
and $F: TQ \rightarrow T^*Q,$ and the control
$u: TQ \rightarrow W (\subset T^*Q)$. This expression (1.1) can not be an equation,
because the left side of (1.1) is defined on the second order
tangent bundle $T^{(2)}Q$, and the right side of (1.1) is defined on the tangent
bundle $TQ$, and $T^{(2)}Q$ and $TQ$ are different spaces. Thus, it is
impossible to define the CL system by using a wrong expression. The
similar wrong appears in the definition 2.4 of \cite{chma04} for the
reduced CL system. In addition, it is worthy of noting that
the use of the above wrong expression (1.1)
has led to the wrong of the method of controlled Lagrangians to judge
the stabilization of mechanical systems,
see Bloch et al. \cite{bllema00, bllema01, blchlema01, blle02},
and this is also a serious problem should to be corrected carefully.\\

(2) The authors didn't consider the phase spaces of CH system and
the reduced CH system, that is, all of CH systems and the reduced
CH systems given in \cite{chbllemawo02,chma04}, have not the spaces
on which these systems are defined, see Definition 3.1 in
\cite{chbllemawo02} and Definition 3.1, 3.3 in \cite{chma04}. Thus,
it is impossible to give the actions of a Lie group
on the phase spaces of CH systems and their momentum maps,
also impossible to determine precisely the
reduced phase spaces of CH systems.\\

(3) The authors didn't consider the change of geometrical
structures of the phase spaces of the CH systems. In fact,
it is not that all of CH systems in \cite{chma04} have same phase space $T^*Q$,
same action of Lie group $G$, and same reduced phase space $T^*Q /G$.
Different structures of geometry determine different CH systems
and their phase spaces.\\

(4) The authors didn't consider the momentum map of the CH system with
symmetry, and hence cannot determine precisely the geometrical
structures of phase spaces of the reduced CH systems.
For example, we consider the cotangent bundle $T^*Q$ of a
smooth manifold $Q$ with a free and proper action of Lie group $G$,
and the Poisson tensor $B$ on $T^*Q$ is determined by the canonical
symplectic form $\omega_0$ on $T^*Q$. Assume that there is an
$\operatorname{Ad}^\ast$-equivariant momentum map $\mathbf{J}: T^*Q
\rightarrow \mathfrak{g}^\ast$ for the symplectic, free and proper
cotangent lifted $G$-action, where $\mathfrak{g}^\ast$ is the dual
of Lie algebra $\mathfrak{g}$ of $G$. For $\mu \in
\mathfrak{g}^\ast$, a regular value of $\mathbf{J}$,
from Abraham and Marsden \cite{abma78}, we know that
the regular point reduced space
$\mathbf{J}^{-1}(\mu)/G_\mu$ and regular orbit reduced space
$\mathbf{J}^{-1}(\mathcal{O}_\mu)/G$ at $\mu$ are different,
and they are not yet the space $T^*Q /G$.
The two reduced spaces are determined by the momentum map
$\mathbf{J}$, where $G_\mu$ is the isotropy subgroup of the
coadjoint $G$-action at the point $\mu$, and $\mathcal{O}_\mu$ is
the orbit of the coadjoint $G$-action through the point $\mu$.
Thus, in the two cases, it is impossible to determine
the reduced CH systems by using the method given in \cite{chma04}.
Moreover, it is also impossible to give precisely the relations of the
reduced controlled Hamiltonian equivalences, if don't consider
the different Lie group actions and momentum maps.\\

To sum up the above statement, we think that there
are a lot of wrong of rigor for the definitions
of CL system, CH system and its reduced CH systems, as well as
CH-equivalence and the reduced CH-equivalence
in Chang et al. \cite{chbllemawo02, chma04},
and we want to correct their work.
It is important to find these wrong, but the more important
is to correct well these wrong. In Marsden et al. \cite{mawazh10},
we have corrected and renewed carefully some of these wrong definitions.\\

In order to deal with the above problems, and give the proper definition
of CH system, and determine uniformly the reduced
CH systems, our idea in Marsden et al. \cite{mawazh10}, is that we first define a CH
system on $T^*Q$ by using the symplectic form, and such system is
called a regular controlled Hamiltonian (RCH) system,
and then regard a Hamiltonian system on $T^*Q$
as a spacial case of an RCH system without external force and
control. Thus, the set of Hamiltonian systems on $T^*Q$ is a subset
of the set of RCH systems on $T^*Q$. On the other hand,
since the symplectic reduced space on a cotangent bundle
is not complete under the Marsden-Weinstein reduction,
and the symplectic reduced system of a Hamiltonian system with symmetry
defined on the cotangent bundle $T^*Q$ may not be a Hamiltonian system
on a cotangent bundle. So,
we can not define directly an RCH system on the cotangent bundle $T^* Q$.
However, from the classification of symplectic
reduced space of the cotangent bundle $T^* Q$, we know that
the regular point symplectic
reduced space $((T^*Q)_\mu, \omega_\mu)$ is symplectically
diffeomorphic to a symplectic fiber bundle over $T^\ast (Q/G_\mu)$
with fiber to be the coadjoint orbit $\mathcal{O}_\mu$. Moreover,
from the regular reduction
diagram, see Ortega and Ratiu \cite{orra04},
we know that the regular orbit reduced space $((T^\ast
Q)_{\mathcal{O}_\mu},\omega_{\mathcal{O}_\mu})$ is symplectically
diffeomorphic to the regular point reduced space $((T^*Q)_\mu,
\omega_\mu)$, and hence is also symplectically diffeomorphic to a
symplectic fiber bundle.
In consequence, if we may define an RCH system on a symplectic fiber bundle,
then it is possible to describe uniformly the RCH system on $T^*Q$ and its
regular reduced RCH systems on the associated reduced spaces, and
we can study regular
reduction theory of the RCH systems with symplectic structures and
symmetries, as an extension of the regular symplectic reduction theory of
Hamiltonian systems under regular controlled Hamiltonian equivalence
conditions. This is why the authors in Marsden et al. \cite{mawazh10}
set up the regular reduction theory of the RCH system on a symplectic fiber
bundle, by using momentum map and the associated reduced symplectic
form and from the viewpoint of completeness of regular symplectic
reduction.\\

Let $(E,M,N,\pi,G)$ be a fiber bundle and $(E, \omega_E)$ be a
symplectic fiber bundle. If for any Hamiltonian function $H: E \rightarrow
\mathbb{R}$, we have a Hamiltonian vector field $X_H$,
which satisfies the Hamilton's equation, that is,
$\mathbf{i}_{X_H}\omega_E=\mathbf{d}H$, then $(E, \omega_E, H )$ is a
Hamiltonian system. Moreover, if considering the external force and
control, we can define a kind of regular controlled Hamiltonian
(RCH) system on the symplectic fiber bundle $E$ as follows.

\begin{defi} (RCH System) An RCH system on $E$ is a 5-tuple
$(E, \omega_E, H, F, W)$, where $(E, \omega_E, H )$ is a
Hamiltonian system, and the function $H: E \rightarrow \mathbb{R}$
is called the Hamiltonian, a fiber-preserving map $F: E\rightarrow
E$ is called the (external) force map, and a fiber sub-manifold
$W$ of $E$ is called the control subset.
\end{defi}
Sometimes, $W$ is also denoted the set of fiber-preserving
maps from $E$ to $W$. When a feedback control law $u:
E\rightarrow W$ is chosen, the 5-tuple $(E, \omega_E, H,
F, u)$ can be denoted a closed-loop dynamical system. In particular, when $Q$
is a smooth manifold, and $T^\ast Q$ its cotangent bundle with a
symplectic form $\omega$ (not necessarily canonical symplectic
form), then $(T^\ast Q, \omega )$ is a symplectic vector bundle. If
we take that $E= T^* Q$, from above definition we can obtain an RCH
system on the cotangent bundle $T^\ast Q$, that is, 5-tuple $(T^\ast
Q, \omega, H, F, W)$. Here for convenience, we assume that all
controls appearing in this paper are the admissible controls.\\

The main contributions in Marsden et al. \cite{mawazh10} are given as
follows. (1) In order to describe uniformly RCH systems defined on a
cotangent bundle and on the regular reduced spaces, we define a kind
of RCH systems on a symplectic fiber bundle by using its symplectic
form; (2) We give regular point and regular orbit reducible RCH
systems by using momentum maps and the associated reduced symplectic
forms, and prove regular point and regular orbit reduction theorems
for the RCH systems, which explain the relationships
between RpCH-equivalence, RoCH-equivalence for the reducible RCH systems
with symmetries and RCH-equivalence for the associated reduced RCH
systems; (3) We prove that rigid body with external force torque,
rigid body with internal rotors and heavy top with internal rotors
are all RCH systems, and as a pair of the regular point reduced RCH
systems, rigid body with internal rotors (or external force torque)
and heavy top with internal rotors are RCH-equivalent; (4) We
describe the RCH system from the viewpoint of port Hamiltonian
system with a symplectic structure, and state the relationship
between RCH-equivalence of RCH systems and equivalence of port
Hamiltonian systems.

\section{The Dynamics of Controlled Hamiltonian System}

In order to describe the dynamics of an RCH system, we have to give a
good expression of the dynamical vector field of the RCH system, by using
the notations of vertical lifted maps of a vector along a fiber, see
Marsden et al. \cite{mawazh10}.\\

At first, for the notations of vertical lifts along fiber,
we need to consider three case:
(1) $\pi: E \rightarrow M$ is a fiber bundle;
(2) $\pi: E \rightarrow M$ is a vector bundle;
(3) $\pi: E \rightarrow M, \; E= T^*Q, \; M=Q, $ is
a cotangent bundle, which is a special vector bundle.
For the case (2) and (3), we can use the standard definition
of the vertical lift operator given in Marsden and Ratiu \cite{mara99}.
But for the case (1), the above operator cannot be used.
This question is found by one of referees who give us that
in a review report of our manuscript.
In order to deal with uniformly
the three cases, we have to give a new definition of vertical
lifted maps of a vector along a fiber, and make it to be
not conflict with that given in Marsden and Ratiu \cite{mara99},
and it is not and cannot be an
extension of the definition of Marsden and Ratiu.\\

It is worthy of noting that there are two aspects in our new definition.
First, for two different points, $a_x,\; b_x $ in the fiber $E_x$, how
define the moving vertical part of a vector in one point $b_x$ to
another point $a_x$; Second, for a fiber-preserving map
$F: E \rightarrow E, $ we know that $a_x$ and $F_x(a_x)$
are the two points in $E_x$,
how define the moving vertical part of a tangent vector in image
point $F_x(a_x)$ to $a_x$. The eventual goal is to give a good
expression of the dynamical vector field of RCH system by using the
notation of vertical lift map of a vector along a fiber. Our
definitions are reasonable and clear, and should be stated
explicitly as follows.\\

For a smooth manifold $E$, its tangent bundle $TE$ is a vector
bundle, and for the fiber bundle $\pi: E \rightarrow M$, we consider
the tangent mapping $T\pi: TE \rightarrow TM$ and its kernel $ker
(T\pi)=\{\rho\in TE| T\pi(\rho)=0\}$, which is a vector subbundle of
$TE$. Denote $VE:= ker(T\pi)$, which is called a vertical bundle
of $E$. Assume that there is a metric on $E$, and we take a
Levi-Civita connection $\mathcal{A}$ on $TE$, and denote by $HE:=
ker(\mathcal{A})$, which is called a horizontal bundle of $E$, such
that $TE= HE \oplus VE. $ For any $x\in M, \; a_x, b_x \in E_x, $
any tangent vector $\rho(b_x)\in T_{b_x}E$ can be split into
horizontal and vertical parts, that is, $\rho(b_x)=
\rho^h(b_x)\oplus \rho^v(b_x)$, where $\rho^h(b_x)\in H_{b_x}E$ and
$\rho^v(b_x)\in V_{b_x}E$. Let $\gamma$ is a geodesic in $E_x$
connecting $a_x$ and $b_x$, and denote by $\rho^v_\gamma(a_x)$ a
tangent vector at $a_x$, which is a parallel displacement of the
vertical vector $\rho^v(b_x)$ along the geodesic $\gamma$ from $b_x$
to $a_x$. Since the angle between two vectors is invariant under a
parallel displacement along a geodesic, then
$T\pi(\rho^v_\gamma(a_x))=0, $ and hence $\rho^v_\gamma(a_x) \in
V_{a_x}E. $ Now, for $a_x, b_x \in E_x $ and tangent vector
$\rho(b_x)\in T_{b_x}E$, we can define the vertical lift map of a
vector along a fiber given by
$$\mbox{vlift}: TE_x \times E_x \rightarrow TE_x; \;\; \mbox{vlift}(\rho(b_x),a_x) = \rho^v_\gamma(a_x). $$
It is easy to check from the basic fact in differential geometry
that this map does not depend on the choice of the geodesic $\gamma$.\\

If $F: E \rightarrow E$ is a fiber-preserving map, for any $x\in M$, we have
that $F_x: E_x \rightarrow E_x$ and $TF_x: TE_x \rightarrow TE_x$,
then for any $a_x \in E_x$ and $\rho\in TE_x$, the vertical lift of
$\rho$ under the action of $F$ along a fiber is defined by
$$(\mbox{vlift}(F_x)\rho)(a_x)=\mbox{vlift}((TF_x\rho)(F_x(a_x)), a_x)= (TF_x\rho)^v_\gamma(a_x), $$
where $\gamma$ is a geodesic in $E_x$ connecting $F_x(a_x)$ and
$a_x$.\\

In particular, when $\pi: E \rightarrow M$ is a vector bundle, for
any $x\in M$, the fiber $E_x$ is a vector space. In this case, we
can choose the geodesic $\gamma$ to be a straight line, and the
vertical vector is invariant under a parallel displacement along a
straight line, that is, $\rho^v_\gamma(a_x)= \rho^v(b_x).$ Moreover,
when $E= T^*Q$, by using the local trivialization of $TT^*Q$, we
have that $TT^*Q\cong TQ \times T^*Q$. Note that $\pi: T^*Q
\rightarrow Q$, and $T\pi: TT^*Q \rightarrow TQ$, then in this case,
for any $\alpha_x, \; \beta_x \in T^*_x Q, \; x\in Q, $ we know that
$(0, \beta_x) \in V_{\beta_x}T^*_x Q, $ and hence we can get that
$$ \mbox{vlift}((0, \beta_x)(\beta_x), \alpha_x) = (0, \beta_x)(\alpha_x), $$
which is consistent with the definition of vertical lift operator along
a fiber given in Marsden and Ratiu \cite{mara99}.\\

For a given RCH System $(T^\ast Q, \omega, H, F, W)$, the dynamical
vector field $X_H$ of the associated Hamiltonian system $(T^\ast Q,
\omega, H) $ satisfies $\mathbf{i}_{X_H}\omega=\mathbf{d}H $.
If considering the external force $F: T^*Q \rightarrow T^*Q, $
which is a fiber-preserving map, by using the above
notations of vertical lift maps of a vector along a fiber, the
change of $X_H$ under the action of $F$ is that
$$\mbox{vlift}(F)X_H(\alpha_x)= \mbox{vlift}((TFX_H)(F(\alpha_x)), \alpha_x)= (TFX_H)^v_\gamma(\alpha_x),$$
where $\alpha_x \in T^*_x Q, \; x\in Q $ and the geodesic $\gamma$ is a straight
line in $T^*_x Q$ connecting $F_x(\alpha_x)$ and $\alpha_x$. In the
same way, when a feedback control law $u: T^\ast Q \rightarrow W,$
which is a fiber-preserving map, is chosen,
the change of $X_H$ under the action of $u$ is that
$$\mbox{vlift}(u)X_H(\alpha_x)= \mbox{vlift}((TuX_H)(F(\alpha_x)), \alpha_x)= (TuX_H)^v_\gamma(\alpha_x).$$
In consequence, we can give an expression of the dynamical vector
field of the RCH system as follows.
\begin{theo}
The dynamical vector field of an RCH system $(T^\ast
Q,\omega,H,F,W)$ with a control law $u$ is the synthetic
of Hamiltonian vector field $X_H$ and its changes under the actions
of the external force $F$ and control $u$, that is,
$$X_{(T^\ast Q,\omega,H,F,u)}(\alpha_x)= X_H(\alpha_x)+ \textnormal{vlift}(F)X_H(\alpha_x)
+ \textnormal{vlift}(u)X_H(\alpha_x),$$ for any $\alpha_x \in T^*_x
Q, \; x\in Q $. For convenience, it is simply written as
\begin{equation}X_{(T^\ast Q,\omega,H,F,u)}
=X_H +\textnormal{vlift}(F) +\textnormal{vlift}(u).
\end{equation}
\end{theo}
We also denote that $\mbox{vlift}(W)=
\bigcup\{\mbox{vlift}(u)X_H | \; u\in W\}$. It is
worthy of noting that in order to deduce and calculate easily, we
always use the simple expressions of the dynamical vector
field $X_{(T^\ast Q,\omega,H,F,u)}$ and the
$R_p$-reduced vector field $X_{((T^\ast
Q)_\mu, \omega_\mu, h_\mu, f_\mu, u_\mu)}$ and the $R_o$-reduced vector
field $X_{((T^\ast Q)_{\mathcal{O}_\mu}, \omega_{\mathcal{O}_\mu},
h_{\mathcal{O}_\mu},f_{\mathcal{O}_\mu},u_{\mathcal{O}_\mu})}$.\\

From the expression (4.1) of the dynamical vector
field of an RCH system, we know that under the actions of the external force $F$
and control $u$, in general, the dynamical vector
field is not Hamiltonian, and hence the RCH system is not
yet a Hamiltonian system. However,
it is a dynamical system closed relative to a
Hamiltonian system, and it can be explored and studied by extending
the methods for external force and control
in the study of Hamiltonian system.

\section{Controlled Hamiltonian Equivalence}

It is worthy of noting that, when an RCH system is given, the force map $F$ is
determined, but the feedback control law $u: T^\ast Q\rightarrow
W $ could be chosen. In order to emphasize explicitly the impact of external force
and control in the study of the RCH systems, by using
the above expression of the dynamical vector field of the RCH system,
we can describe the feedback
control law how to modify the structure of an RCH system, and the controlled
Hamiltonian matching conditions and RCH-equivalence are induced as
follows.
\begin{defi}
(RCH-equivalence) Suppose that we have two RCH systems $(T^\ast
Q_i,\omega_i,H_i,F_i,W_i),$ $ i= 1,2,$ we say them to be
RCH-equivalent, or simply, $(T^\ast
Q_1,\omega_1,H_1,F_1,W_1)\stackrel{RCH}{\sim}\\ (T^\ast
Q_2,\omega_2,H_2,F_2,W_2)$, if there exists a
diffeomorphism $\varphi: Q_1\rightarrow Q_2$, such that the
following controlled Hamiltonian matching conditions hold:

\noindent {\bf RCH-1:} The cotangent lifted map of $\varphi$, that
is, $\varphi^\ast= T^\ast \varphi:T^\ast Q_2\rightarrow T^\ast Q_1$
is symplectic, and $W_1=\varphi^\ast (W_2).$

\noindent {\bf RCH-2:} $Im[X_{H_1}+ \textnormal{vlift}(F_1)-
T\varphi^\ast (X_{H_2})-\textnormal{vlift}(\varphi^\ast
F_2\varphi_\ast)]\subset \textnormal{vlift}(W_1)$, where
the map $\varphi_\ast=(\varphi^{-1})^\ast: T^\ast Q_1\rightarrow
T^\ast Q_2$, and $T\varphi^\ast: TT^\ast Q_2\rightarrow TT^\ast
Q_1$, and $Im$ means the pointwise image of the map in brackets.
\end{defi}

It is worthy of noting that (i) we define directly the RCH
system $(T^\ast Q, \omega, H, F, W)$ with a space $T^*Q$; (ii) we
use a diffeomorphism $\varphi: Q_1\rightarrow Q_2$ to describe that
two different RCH systems $(T^\ast Q_i,\omega_i,H_i,F_i,W_i),\; i=
1,2,$ are RCH-equivalent. These are very important for a rigorous
definition of CH-equivalence, by comparing the definitions in
Chang et al. \cite{chbllemawo02, chma04}.
On the other hand, our RCH system is defined by using the
symplectic structure on the cotangent bundle of a configuration
manifold, we have to keep with the symplectic structure when we define
the RCH-equivalence, that is, the induced equivalent map $\varphi^*$
is symplectic on the cotangent bundle.
Moreover, the following theorem explains the
significance of the above RCH-equivalence relation,
its proof is given in Marsden et al. \cite{mawazh10}.
\begin{theo}
Suppose that two RCH systems $(T^\ast Q_i,\omega_i,H_i,F_i,W_i)$,
$i=1,2,$ are RCH-equivalent, then there exist two control laws $u_i:
T^\ast Q_i \rightarrow W_i, \; i=1,2, $ such that the two
closed-loop dynamical systems produce the same equations of motion, that is,
$X_{(T^\ast Q_1,\omega_1,H_1,F_1,u_1)}\cdot \varphi^\ast
=T(\varphi^\ast) X_{(T^\ast Q_2,\omega_2,H_2,F_2,u_2)}$, where the
map $T(\varphi^\ast):TT^\ast Q_2\rightarrow TT^\ast Q_1$ is the
tangent map of $\varphi^\ast$. Moreover, the explicit relation
between the two control laws $u_i, i=1,2$ is given by
\begin{equation}\textnormal{vlift}(u_1) -\textnormal{vlift}(\varphi^\ast
u_2\varphi_\ast)=-X_{H_1}-\textnormal{vlift}(F_1)
+T\varphi^\ast (X_{H_2})+\textnormal{vlift}(\varphi^\ast F_2
\varphi_\ast)\label{3.2}\end{equation}
\end{theo}

\section{Regular Reducible Controlled Hamiltonian System}

We know that when the external force and control of an RCH
system $(T^*Q,\omega,H,F,W)$ are both zeros, in this case the RCH system
is just a Hamiltonian system $(T^*Q,\omega,H)$.
Thus, we can regard a Hamiltonian system on $T^*Q$
as a spacial case of an RCH system without external force and
control. In consequence, the set of Hamiltonian systems with symmetries
on $T^*Q$ is a subset of the set of RCH systems with symmetries on $T^*Q$.
If we first admit the regular symplectic reduction of a Hamiltonian system
with symmetry, then we may study the regular
reduction of an RCH system with symmetry,
as an extension of regular symplectic reduction of
a Hamiltonian system under regular controlled Hamiltonian equivalence
conditions. In order to do these, in this section
we first give the regular point reducible RCH system
and the regular orbit reducible RCH system,
by using the Marsden-Weinstein reduction and
regular orbit reduction for a Hamiltonian system, respectively.\\

At first, we consider the regular point reducible RCH system.
Let $Q$ be a smooth manifold and $T^\ast Q$ its cotangent bundle
with the symplectic form $\omega$. Let $\Phi: G\times Q\rightarrow
Q$ be a smooth left action of a Lie group $G$ on $Q$, which is free
and proper. Then the cotangent lifted left action $\Phi^{T^\ast}:
G\times T^\ast Q\rightarrow T^\ast Q$ is also free and
proper. Assume that the action is symplectic and admits an
$\operatorname{Ad}^\ast$-equivariant momentum map $\mathbf{J}:T^\ast
Q\rightarrow \mathfrak{g}^\ast$, where $\mathfrak{g}$ is the Lie
algebra of $G$ and $\mathfrak{g}^\ast$ is the dual of
$\mathfrak{g}$. Let $\mu\in\mathfrak{g}^\ast$ be a regular value of
$\mathbf{J}$ and denote by $G_\mu$ the isotropy subgroup of the
coadjoint $G$-action at the point $\mu\in\mathfrak{g}^\ast$, which
is defined by $G_\mu=\{g\in G|\operatorname{Ad}_g^\ast \mu=\mu \}$.
Since $G_\mu (\subset G)$ acts freely and properly on $Q$ and on
$T^\ast Q$, then $Q_\mu=Q/G_\mu$ is a smooth manifold and that the
canonical projection $\rho_\mu:Q\rightarrow Q_\mu$ is a surjective
submersion. It follows that $G_\mu$ acts also freely and properly on
$\mathbf{J}^{-1}(\mu)$, so that the space $(T^\ast
Q)_\mu=\mathbf{J}^{-1}(\mu)/G_\mu$ is a symplectic manifold with the
symplectic form $\omega_\mu$ uniquely characterized by the relation
\begin{equation}\pi_\mu^\ast \omega_\mu=i_\mu^\ast
\omega. \label{6.1}\end{equation} The map
$i_\mu:\mathbf{J}^{-1}(\mu)\rightarrow T^\ast Q$ is the inclusion
and $\pi_\mu:\mathbf{J}^{-1}(\mu)\rightarrow (T^\ast Q)_\mu$ is the
projection. The pair $((T^\ast Q)_\mu,\omega_\mu)$ is called
Marsden-Weinstein reduced space of $(T^\ast Q,\omega)$ at $\mu$.\\

Assume that $H: T^\ast Q\rightarrow \mathbb{R}$ is a $G$-invariant
Hamiltonian, the flow $F_t$ of the Hamiltonian vector field $X_H$
leaves the connected components of $\mathbf{J}^{-1}(\mu)$ invariant
and commutes with the $G$-action, so it induces a flow $f_t^\mu$ on
$(T^\ast Q)_\mu$, defined by $f_t^\mu\cdot \pi_\mu=\pi_\mu \cdot
F_t\cdot i_\mu$, and the vector field $X_{h_\mu}$ generated by the
flow $f_t^\mu$ on $((T^\ast Q)_\mu,\omega_\mu)$ is Hamiltonian with
the associated regular point reduced Hamiltonian function
$h_\mu:(T^\ast Q)_\mu\rightarrow \mathbb{R}$ defined by
$h_\mu\cdot\pi_\mu=H\cdot i_\mu$, and the Hamiltonian vector fields
$X_H$ and $X_{h_\mu}$ are $\pi_\mu$-related. On the other hand, from
Marsden et al. \cite{mawazh10}, we know that the regular point
reduced space $((T^*Q)_\mu, \omega_\mu)$ is symplectically
diffeomorphic to a symplectic fiber bundle. Thus, we can introduce a
kind of the regular point reducible RCH system as follows.
\begin{defi}
(Regular Point Reducible RCH System) A 6-tuple $(T^\ast Q, G,
\omega, H, F, W)$, where the Hamiltonian $H:T^\ast Q\rightarrow
\mathbb{R}$, the fiber-preserving map $F:T^\ast Q\rightarrow T^\ast
Q$ and the fiber submanifold $W$ of\; $T^\ast Q$ are all
$G$-invariant, is called a regular point reducible RCH system, if
there exists a point $\mu\in\mathfrak{g}^\ast$, which is a regular
value of the momentum map $\mathbf{J}$, such that the regular point
reduced system, that is, the 5-tuple $((T^\ast Q)_\mu,
\omega_\mu,h_\mu,f_\mu,W_\mu)$, where $(T^\ast
Q)_\mu=\mathbf{J}^{-1}(\mu)/G_\mu$, $\pi_\mu^\ast
\omega_\mu=i_\mu^\ast\omega$, $h_\mu\cdot \pi_\mu=H\cdot i_\mu$,
$F(\mathbf{J}^{-1}(\mu))\subset \mathbf{J}^{-1}(\mu) $, $f_\mu\cdot
\pi_\mu=\pi_\mu \cdot F\cdot i_\mu$, $W \cap
\mathbf{J}^{-1}(\mu)\neq \emptyset $, $W_\mu=\pi_\mu(W\cap
\mathbf{J}^{-1}(\mu))$, is an RCH system, which is simply written as
$R_p$-reduced RCH system. Where $((T^\ast Q)_\mu,\omega_\mu)$ is the
$R_p$-reduced space, the function $h_\mu:(T^\ast Q)_\mu\rightarrow
\mathbb{R}$ is called the $R_p$-reduced Hamiltonian, the fiber-preserving
map $f_\mu:(T^\ast Q)_\mu\rightarrow (T^\ast Q)_\mu$ is called the
$R_p$-reduced (external) force map, $W_\mu$ is a fiber submanifold of
\;$(T^\ast Q)_\mu$ and is called the $R_p$-reduced control subset.
\end{defi}

It is worthy of noting that for the regular point reducible RCH system
$(T^\ast Q,G,\omega,H,F,W)$, the $G$-invariant external force map
$F: T^*Q \rightarrow T^*Q $ has to satisfy the conditions
$F(\mathbf{J}^{-1}(\mu))\subset \mathbf{J}^{-1}(\mu), $ and
$f_\mu\cdot \pi_\mu=\pi_\mu \cdot F\cdot i_\mu, $ such that we can
define the $R_p$-reduced external force map $f_\mu:(T^\ast
Q)_\mu\rightarrow (T^\ast Q)_\mu. $ The condition $W \cap
\mathbf{J}^{-1}(\mu)\neq \emptyset $ in above definition makes that
the $G$-invariant control subset $W\cap \mathbf{J}^{-1}(\mu)$ can be
reduced and the $R_p$-reduced control subset is $W_\mu= \pi_\mu(W\cap
\mathbf{J}^{-1}(\mu))$. If the control subset cannot be reduced, we
cannot get the $R_p$-reduced RCH system. The study of RCH system
which is not regular point reducible is beyond the limits in this
paper, it may be a topic in future study.\\

Next, we consider the regular orbit reducible RCH system.
For the cotangent lifted left action
$\Phi^{T^\ast}:G\times T^\ast Q\rightarrow T^\ast Q$, which is
symplectic, free and proper, assume that the action admits an
$\operatorname{Ad}^\ast$-equivariant momentum map $\mathbf{J}:T^\ast
Q\rightarrow \mathfrak{g}^\ast$. Let $\mu\in \mathfrak{g}^\ast$ be a
regular value of the momentum map $\mathbf{J}$ and
$\mathcal{O}_\mu=G\cdot \mu\subset \mathfrak{g}^\ast$ be the
$G$-orbit of the coadjoint $G$-action through the point $\mu$. Since
$G$ acts freely, properly and symplectically on $T^\ast Q$, then the
quotient space $(T^\ast Q)_{\mathcal{O}_\mu}=
\mathbf{J}^{-1}(\mathcal{O}_\mu)/G$ is a regular quotient symplectic
manifold with the symplectic form $\omega_{\mathcal{O}_\mu}$
uniquely characterized by the relation
\begin{equation}i_{\mathcal{O}_\mu}^\ast \omega=\pi_{\mathcal{O}_{\mu}}^\ast
\omega_{\mathcal{O}
_\mu}+\mathbf{J}_{\mathcal{O}_\mu}^\ast\omega_{\mathcal{O}_\mu}^+,
\label{6.2}\end{equation} where $\mathbf{J}_{\mathcal{O}_\mu}$ is
the restriction of the momentum map $\mathbf{J}$ to
$\mathbf{J}^{-1}(\mathcal{O}_\mu)$, that is,
$\mathbf{J}_{\mathcal{O}_\mu}=\mathbf{J}\cdot i_{\mathcal{O}_\mu}$
and $\omega_{\mathcal{O}_\mu}^+$ is the $+$-symplectic structure on
the orbit $\mathcal{O}_\mu$ given by
\begin{equation}\omega_{\mathcal{O}_\mu}^
+(\nu)(\xi_{\mathfrak{g}^\ast}(\nu),\eta_{\mathfrak{g}^\ast}(\nu))
=<\nu,[\xi,\eta]>,\;\; \forall\;\nu\in\mathcal{O}_\mu, \;
\xi,\eta\in \mathfrak{g}. \label{6.3}\end{equation} The maps
$i_{\mathcal{O}_\mu}:\mathbf{J}^{-1}(\mathcal{O}_\mu)\rightarrow
T^\ast Q$ and
$\pi_{\mathcal{O}_\mu}:\mathbf{J}^{-1}(\mathcal{O}_\mu)\rightarrow
(T^\ast Q)_{\mathcal{O}_\mu}$ are natural injection and the
projection, respectively. The pair $((T^\ast
Q)_{\mathcal{O}_\mu},\omega_{\mathcal{O}_\mu})$ is called the
symplectic orbit reduced space of $(T^\ast Q,\omega)$ at $\mu$.\\

Assume that $H:T^\ast Q\rightarrow \mathbb{R}$ is a $G$-invariant
Hamiltonian, the flow $F_t$ of the Hamiltonian vector field $X_H$
leaves the connected components of
$\mathbf{J}^{-1}(\mathcal{O}_\mu)$ invariant and commutes with the
$G$-action, so it induces a flow $f_t^{\mathcal{O}_\mu}$ on $(T^\ast
Q)_{\mathcal{O}_\mu}$, defined by $f_t^{\mathcal{O}_\mu}\cdot
\pi_{\mathcal{O}_\mu}=\pi_{\mathcal{O}_\mu} \cdot F_t\cdot
i_{\mathcal{O}_\mu}$, and the vector field $X_{h_{\mathcal{O}_\mu}}$
generated by the flow $f_t^{\mathcal{O}_\mu}$ on $((T^\ast
Q)_{\mathcal{O}_\mu},\omega_{\mathcal{O}_\mu})$ is Hamiltonian with
the associated regular orbit reduced Hamiltonian function
$h_{\mathcal{O}_\mu}:(T^\ast Q)_{\mathcal{O}_\mu}\rightarrow
\mathbb{R}$ defined by $h_{\mathcal{O}_\mu}\cdot
\pi_{\mathcal{O}_\mu}= H\cdot i_{\mathcal{O}_\mu}$, and the
Hamiltonian vector fields $X_H$ and $X_{h_{\mathcal{O}_\mu}}$ are
$\pi_{\mathcal{O}_\mu}$-related. In general case, we
maybe thought that the structure of the symplectic orbit reduced
space $((T^\ast Q)_{\mathcal{O}_\mu},\omega_{\mathcal{O}_\mu})$ is
more complex than that of the symplectic point reduced space
$((T^\ast Q)_\mu,\omega_\mu)$, but, from the regular reduction
diagram, see Ortega and Ratiu \cite{orra04},
we know that the regular orbit reduced space $((T^\ast
Q)_{\mathcal{O}_\mu},\omega_{\mathcal{O}_\mu})$ is symplectically
diffeomorphic to the regular point reduced space $((T^*Q)_\mu,
\omega_\mu)$, and hence is also symplectically diffeomorphic to a
symplectic fiber bundle. Thus, we can introduce a kind of the
regular orbit reducible RCH systems as follows.
\begin{defi}
(Regular Orbit Reducible RCH System) A 6-tuple $(T^\ast Q, G,
\omega,H,F,W)$, where the Hamiltonian $H: T^\ast Q\rightarrow
\mathbb{R}$, the fiber-preserving map $F: T^\ast Q\rightarrow T^\ast
Q$ and the fiber submanifold $W$ of $T^\ast Q$ are all
$G$-invariant, is called a regular orbit reducible RCH system, if
there exists an orbit $\mathcal{O}_\mu, \; \mu\in\mathfrak{g}^\ast$,
where $\mu$ is a regular value of the momentum map $\mathbf{J}$,
such that the regular orbit reduced system, that is, the 5-tuple
$((T^\ast
Q)_{\mathcal{O}_\mu},\omega_{\mathcal{O}_\mu},h_{\mathcal{O}_\mu},f_{\mathcal{O}_\mu},
W_{\mathcal{O}_\mu})$, where $(T^\ast
Q)_{\mathcal{O}_\mu}=\mathbf{J}^{-1}(\mathcal{O}_\mu)/G$,
$\pi_{\mathcal{O}_\mu}^\ast \omega_{\mathcal{O}_\mu}
=i_{\mathcal{O}_\mu}^\ast\omega-\mathbf{J}_{\mathcal{O}_\mu}^\ast\omega_{\mathcal{O}_\mu}^+$,
$h_{\mathcal{O}_\mu}\cdot \pi_{\mathcal{O}_\mu} =H\cdot
i_{\mathcal{O}_\mu}$, $F(\mathbf{J}^{-1}(\mathcal{O}_\mu))\subset
\mathbf{J}^{-1}(\mathcal{O}_\mu)$, $f_{\mathcal{O}_\mu}\cdot
\pi_{\mathcal{O}_\mu}=\pi_{\mathcal{O}_\mu}\cdot F\cdot
i_{\mathcal{O}_\mu}$, and $W \cap
\mathbf{J}^{-1}(\mathcal{O}_\mu)\neq \emptyset $,
$W_{\mathcal{O}_\mu}=\pi_{\mathcal{O}_\mu}(W \cap
\mathbf{J}^{-1}(\mathcal{O}_\mu))$, is an RCH system, which is simply
written as the $R_o$-reduced RCH system. Where $((T^\ast
Q)_{\mathcal{O}_\mu},\omega_{\mathcal{O}_\mu})$ is the $R_o$-reduced
space, the function $h_{\mathcal{O}_\mu}:(T^\ast
Q)_{\mathcal{O}_\mu}\rightarrow \mathbb{R}$ is called the $R_o$-reduced
Hamiltonian, the fiber-preserving map $f_{\mathcal{O}_\mu}:(T^\ast
Q)_{\mathcal{O}_\mu} \rightarrow (T^\ast Q)_{\mathcal{O}_\mu}$ is
called the $R_o$-reduced (external) force map, $W_{\mathcal{O}_\mu}$ is a
fiber submanifold of $(T^\ast Q)_{\mathcal{O}_\mu}$, and is called
the $R_o$-reduced control subset.
\end{defi}

It is worthy of noting that for the regular orbit reducible RCH system
$(T^\ast Q,G,\omega,H,F,W)$, the $G$-invariant external force map
$F: T^*Q \rightarrow T^*Q $ has to satisfy the conditions
$F(\mathbf{J}^{-1}(\mathcal{O}_\mu))\subset
\mathbf{J}^{-1}(\mathcal{O}_\mu), $ and $f_{\mathcal{O}_\mu}\cdot
\pi_{\mathcal{O}_\mu}=\pi_{\mathcal{O}_\mu}\cdot F\cdot
i_{\mathcal{O}_\mu}$, such that we can define the reduced external
force map $f_{\mathcal{O}_\mu}:(T^\ast Q)_{\mathcal{O}_\mu}
\rightarrow (T^\ast Q)_{\mathcal{O}_\mu}. $ The condition $W \cap
\mathbf{J}^{-1}(\mathcal{O}_\mu)\neq \emptyset $ in above definition
makes that the  $G$-invariant control subset $W \cap
\mathbf{J}^{-1}(\mathcal{O}_\mu)$ can be reduced and the reduced
control subset is $W_{\mathcal{O}_\mu}= \pi_{\mathcal{O}_\mu}(W \cap
\mathbf{J}^{-1}(\mathcal{O}_\mu))$. If the control subset cannot be
reduced, we cannot get the $R_o$-reduced RCH system. The study of
RCH system which is not regular orbit reducible is beyond the limits
in this paper, it may be a topic in future study.

\section{Regular Point Reduction of the RCH System}

In the following we consider the RCH system with symmetry and momentum map,
and give the RpCH-equivalence for the regular point reducible RCH system,
and prove the regular point reduction theorem.
Denote by $X_{(T^\ast Q,G,\omega,H,F,u)}$ the dynamical vector field of the
regular point reducible RCH system $(T^\ast Q,G,\omega,
H,F,W)$ with a control law $u$, and assume that it can be expressed by
\begin{equation}X_{(T^\ast Q,G,\omega,H,F,u)}
=X_H+\textnormal{vlift}(F)+\textnormal{vlift}(u).\label{3.2}\end{equation}
Moreover, by using the above expression,
for the regular point reducible RCH system we can also
introduce the regular point reducible controlled Hamiltonian
equivalence (RpCH-equivalence) as follows.

\begin{defi}(RpCH-equivalence)
Suppose that we have two regular point reducible RCH systems
$(T^\ast Q_i, G_i,\omega_i,H_i, F_i, W_i),\; i=1,2$, we
say them to be RpCH-equivalent, or simply,\\ $(T^\ast Q_1,
G_1,\omega_1,H_1,F_1,W_1)\stackrel{RpCH}{\sim}(T^\ast
Q_2,G_2,\omega_2,H_2,F_2,W_2)$, if there exists a
diffeomorphism $\varphi:Q_1\rightarrow Q_2$ such that the following
controlled Hamiltonian matching conditions hold:

\noindent {\bf RpCH-1:} The cotangent lifted map
$\varphi^\ast:T^\ast Q_2\rightarrow T^\ast Q_1$ is symplectic.

\noindent {\bf RpCH-2:} For $\mu_i\in \mathfrak{g}^\ast_i $, the
regular reducible points of RCH systems $(T^\ast Q_i, G_i,\omega_i,
H_i, F_i, W_i),\\ i=1,2$, the map
$\varphi_\mu^\ast=i_{\mu_1}^{-1}\cdot\varphi^\ast\cdot i_{\mu_2}:
\mathbf{J}_2^{-1}(\mu_2)\rightarrow \mathbf{J}_1^{-1}(\mu_1)$ is
$(G_{2\mu_2},G_{1\mu_1})$-equivariant,
and $W_1\cap \mathbf{J}_1^{-1}(\mu_1) = \varphi_\mu^\ast (W_2\cap
\mathbf{J}_2^{-1}(\mu_2))$, where $\mu=(\mu_1,
\mu_2)$, and denote by $i_{\mu_1}^{-1}(S)$ the pre-image of a subset
$S\subset T^\ast Q_1$ for the map
$i_{\mu_1}:\mathbf{J}_1^{-1}(\mu_1)\rightarrow T^\ast Q_1$.

\noindent {\bf RpCH-3:} $Im[X_{H_1}+ \textnormal{vlift}(F_1)-
T\varphi^\ast (X_{H_2})-\textnormal{vlift}(\varphi^\ast
F_2\varphi_\ast)]\subset\textnormal{vlift}(W_1)$.
\end{defi}

It is worthy of noting that for the regular point reducible RCH
system, the induced equivalent map $\varphi^*$ not only keeps the
symplectic structure, but also keeps the equivariance of $G$-action
at the regular point. If an $R_p$-reduced feedback control law $u_\mu:(T^\ast
Q)_\mu\rightarrow W_\mu$ is chosen, the $R_p$-reduced RCH
system $((T^\ast Q)_\mu, \omega_\mu, h_\mu, f_\mu, u_\mu)$ is a
closed-loop regular dynamic system with a control law $u_\mu$.
Assume that its vector field $X_{((T^\ast Q)_\mu, \omega_\mu, h_\mu,
f_\mu, u_\mu)}$ can be expressed by
\begin{equation}X_{((T^\ast Q)_\mu, \omega_\mu, h_\mu, f_\mu, u_\mu)}
=X_{h_\mu}+\textnormal{vlift}(f_\mu)+\textnormal{vlift}(u_\mu),
\label{3.3}\end{equation}
where $X_{h_\mu}$ is the dynamical vector field of
the $R_p$-reduced Hamiltonian $h_\mu$, $\textnormal{vlift}(f_\mu)=
\textnormal{vlift}(f_\mu)X_{h_\mu}$, $\textnormal{vlift}(u_\mu)=
\textnormal{vlift}(u_\mu)X_{h_\mu}$, and satisfies the condition
\begin{equation}X_{((T^\ast Q)_\mu, \omega_\mu, h_\mu, f_\mu,
u_\mu)}\cdot \pi_\mu=T\pi_\mu\cdot X_{(T^\ast
Q,G,\omega,H,F,u)}\cdot i_\mu. \label{3.4}\end{equation}
Then we can
obtain the following regular point reduction theorem for the RCH
system, which explains the relationship between the RpCH-equivalence
for the regular point reducible RCH systems with symmetries and the
RCH-equivalence for the associated $R_p$-reduced RCH systems, its
proof is given in Marsden et al. \cite{mawazh10}.
\begin{theo}
Two regular point reducible RCH systems $(T^\ast Q_i, G_i, \omega_i,
H_i, F_i,W_i)$, $i=1,2,$ are RpCH-equivalent if and only
if the associated $R_p$-reduced RCH systems $((T^\ast
Q_i)_{\mu_i},\omega_{i\mu_i},h_{i\mu_i},f_{i\mu_i},\\
W_{i\mu_i}),$ $i=1,2,$ are RCH-equivalent.
\end{theo}
This theorem can
be regarded as an extension of the regular point reduction theorem
of Hamiltonian system under the regular controlled Hamiltonian
equivalence conditions.

\section{Regular Orbit Reduction of the RCH System}

The orbit reduction of a Hamiltonian system is an alternative
approach to symplectic reduction given by Marle \cite{ma76}
and Kazhdan, Kostant and Sternberg \cite{kakost78}, which is different from the
Marsden-Weinstein reduction. We note that the regular reduced
symplectic spaces $((T^\ast Q)_{\mathcal{O}_\mu},\omega_{\mathcal{O}_\mu})$
and $((T^\ast Q)_\mu,\omega_\mu), $ of the regular orbit
reduced Hamiltonian system and the regular point reduced Hamiltonian system,
are different, and the symplectic forms on the reduced spaces,
given by (6.2) for the regular orbit
reduced Hamiltonian system and given by (6.1) for the regular point
reduced Hamiltonian system, are also different. Thus, the assumption conditions
for the regular orbit reduction case are not same as
that for the regular point reduction case.\\

In the following we consider the RCH system with symmetry and momentum map,
and give the RoCH-equivalence for the regular orbit reducible RCH systems,
and prove the regular orbit reduction theorem.
Denote by $X_{(T^\ast Q,G,\omega,H,F,u)}$ the dynamical vector field of the
regular orbit reducible RCH system $(T^\ast Q, G,\omega, H,F,W)$
with a control law $u$, and assume that it can be expressed by
\begin{equation}X_{(T^\ast Q,G,\omega,H,F,u)}
=X_H+\textnormal{vlift}(F)+\textnormal{vlift}(u).\label{5.3}
\end{equation}
Moreover, by using the above expression,
for the regular orbit reducible RCH system we can also
introduce the regular orbit reducible controlled Hamiltonian
equivalence (RoCH-equivalence) as follows.
\begin{defi}
(RoCH-equivalence) Suppose that we have two regular orbit reducible
RCH systems $(T^\ast Q_i, G_i, \omega_i, H_i, F_i, W_i)$, $i=1,2$,
we say them to be RoCH-equivalent, or simply,\\ $(T^\ast Q_1, G_1,
\omega_1, H_1, F_1, W_1)\stackrel{RoCH}{\sim}(T^\ast Q_2, G_2,
\omega_2, H_2, F_2, W_2)$, if there exists a diffeomorphism
$\varphi:Q_1\rightarrow Q_2$ such that the following controlled Hamiltonian
matching conditions hold:

\noindent {\bf RoCH-1:} The cotangent lift map $\varphi^\ast:
T^\ast Q_2\rightarrow T^\ast Q_1$ is symplectic.

\noindent {\bf RoCH-2:} For $\mathcal{O}_{\mu_i},\; \mu_i\in
\mathfrak{g}^\ast_i$, the regular reducible orbits of RCH systems
$(T^\ast Q_i, G_i, \omega_i, H_i, F_i, W_i)$, $i=1,2$, the map
$\varphi^\ast_{\mathcal{O}_\mu}=i_{\mathcal{O}_{\mu_1}}^{-1}\cdot\varphi^\ast\cdot
i_{\mathcal{O}_{\mu_2}}:\mathbf{J}_2^{-1}(\mathcal{O}_{\mu_2})\rightarrow
\mathbf{J}_1^{-1}(\mathcal{O}_{\mu_1})$ is $(G_2,G_1)$-equivariant,
 and $W_1\cap \mathbf{J}_1^{-1}(\mathcal{O}_{\mu_1})=\varphi_{\mathcal{O}_\mu}^\ast
(W_2\cap \mathbf{J}_2^{-1}(\mathcal{O}_{\mu_2}))$, and
$\mathbf{J}_{2\mathcal{O}_{\mu_2}}^\ast
\omega_{2\mathcal{O}_{\mu_2}}^{+}=(\varphi_{\mathcal{O}_\mu}^\ast)^\ast
\cdot\mathbf{J}_{1\mathcal{O}_{\mu_1}}^\ast\omega_{1\mathcal{O}_{\mu_1}}^{+},$
where $\mu=(\mu_1, \mu_2)$, and denote by
$i_{\mathcal{O}_{\mu_1}}^{-1}(S)$ the pre-image of a subset
$S\subset T^\ast Q_1$ for the map
$i_{\mathcal{O}_{\mu_1}}:\mathbf{J}_1^{-1}(\mathcal{O}_{\mu_1})\rightarrow
T^\ast Q_1$.

\noindent {\bf RoCH-3:}
$Im[X_{H_1}+\textnormal{vlift}(F_1)-T\varphi^\ast (X_{H_2})
-\textnormal{vlift}(\varphi^\ast F_2\varphi_\ast)]\subset
\textnormal{vlift}(W_1).$
\end{defi}

It is worthy of noting that for the regular orbit reducible RCH
system, the induced equivalent map $\varphi^*$ not only keeps the
symplectic structure and the restriction of the $(+)$-symplectic
structure on the regular orbit to
$\mathbf{J}^{-1}(\mathcal{O}_\mu)$, but also keeps the equivariance
of $G$-action on the regular orbit. If an $R_o$-reduced feedback control law
$u_{\mathcal{O}_\mu}:(T^\ast Q)_{\mathcal{O}_\mu}\rightarrow
W_{\mathcal{O}_\mu}$ is chosen, the $R_o$-reduced RCH system
$((T^\ast Q)_{\mathcal{O}_\mu},
\omega_{\mathcal{O}_\mu},h_{\mathcal{O}_\mu},f_{\mathcal{O}_\mu},u_{\mathcal{O}_\mu})$
is a closed-loop regular dynamic system with a control law
$u_{\mathcal{O}_\mu}$. Assume that its vector field $X_{((T^\ast
Q)_{\mathcal{O}_\mu}, \omega_{\mathcal{O}_\mu},
h_{\mathcal{O}_\mu},f_{\mathcal{O}_\mu},u_{\mathcal{O}_\mu})}$ can
be expressed by
\begin{equation}X_{((T^\ast Q)_{\mathcal{O}_\mu},
\omega_{\mathcal{O}_\mu},h_{\mathcal{O}_\mu},f_{\mathcal{O}_\mu},u_{\mathcal{O}_\mu})}=
X_{h_{\mathcal{O}_\mu}}+\textnormal{vlift}(f_{\mathcal{O}_\mu})
+\textnormal{vlift}(u_{\mathcal{O}_\mu}), \label{5.4}\end{equation}
where $X_{h_{\mathcal{O}_\mu}}$ is the dynamical vector field of
the $R_o$-reduced Hamiltonian $h_{\mathcal{O}_\mu}$, and $\textnormal{vlift}(f_{\mathcal{O}_\mu})=
\textnormal{vlift}(f_{\mathcal{O}_\mu})X_{h_{\mathcal{O}_\mu}}$,
$\textnormal{vlift}(u_{\mathcal{O}_\mu})=
\textnormal{vlift}(u_{\mathcal{O}_\mu})X_{h_{\mathcal{O}_\mu}}$, and
satisfies the condition
\begin{equation}X_{((T^\ast Q)_{\mathcal{O}_\mu},
\omega_{\mathcal{O}_\mu},h_{\mathcal{O}_\mu},f_{\mathcal{O}_\mu},u_{\mathcal{O}_\mu})}\cdot
\pi_{\mathcal{O}_\mu} =T\pi_{\mathcal{O}_\mu} \cdot X_{(T^\ast
Q,G,\omega,H,F,u)}\cdot
i_{\mathcal{O}_\mu}.\label{5.5}\end{equation}
Then we can obtain the
following regular orbit reduction theorem for the RCH system, which
explains the relationship between the RoCH-equivalence for regular
orbit reducible RCH systems with symmetries and the RCH-equivalence
for the associated $R_o$-reduced RCH systems,
its proof is given in Marsden et al. \cite{mawazh10}.
\begin{theo}
If two regular orbit reducible RCH systems $(T^\ast Q_i, G_i,
\omega_i, H_i, F_i,W_i)$, $i=1,2,$ are RoCH-equivalent, then their
associated $R_o$-reduced RCH systems $((T^\ast
Q)_{\mathcal{O}_{\mu_i}}, \omega_{i\mathcal{O}_{\mu_i}},
h_{i\mathcal{O}_{\mu_i}}, f_{i\mathcal{O}_{\mu_i}},
W_{i\mathcal{O}_{\mu_i}})$, $i=1,2,$ must be RCH-equivalent.
Conversely, if $R_o$-reduced RCH systems $((T^\ast
Q)_{\mathcal{O}_{\mu_i}}, \omega_{i\mathcal{O}_{\mu_i}},
h_{i\mathcal{O}_{\mu_i}},\\ f_{i\mathcal{O}_{\mu_i}},
W_{i\mathcal{O}_{\mu_i}})$, $i=1,2,$ are RCH-equivalent and the
induced map
$\varphi^\ast_{\mathcal{O}_\mu}:\mathbf{J}_2^{-1}(\mathcal{O}_{\mu_2})\rightarrow
\mathbf{J}_1^{-1}(\mathcal{O}_{\mu_1})$, such that
$\mathbf{J}_{2\mathcal{O}_{\mu_2}}^\ast
\omega_{2\mathcal{O}_{\mu_2}}^{+}=(\varphi_{\mathcal{O}_\mu}^\ast)^\ast
\cdot\mathbf{J}_{1\mathcal{O}_{\mu_1}}^\ast\omega_{1\mathcal{O}_{\mu_1}}^{+},$
then the regular orbit reducible RCH systems $(T^\ast Q_i, G_i,
\omega_i, H_i,\\  F_i, W_i)$, $i=1,2,$ are RoCH-equivalent.
\end{theo}
This theorem can be
regarded as an extension of the regular orbit reduction theorem of
Hamiltonian system under the regular controlled Hamiltonian equivalence
conditions.

\begin{rema}
If $(T^\ast Q, \omega)$ is a connected symplectic manifold, and
$\mathbf{J}:T^\ast Q\rightarrow \mathfrak{g}^\ast$ is a
non-equivariant momentum map with a non-equivariance group
one-cocycle $\sigma: G\rightarrow \mathfrak{g}^\ast$, which is
defined by $\sigma(g):=\mathbf{J}(g\cdot
z)-\operatorname{Ad}^\ast_{g^{-1}}\mathbf{J}(z)$, where $g\in G$ and
$z\in T^\ast Q$. Then we know that $\sigma$ produces a new affine
action $\Theta: G\times \mathfrak{g}^\ast \rightarrow
\mathfrak{g}^\ast $ defined by
$\Theta(g,\mu):=\operatorname{Ad}^\ast_{g^{-1}}\mu + \sigma(g)$,
where $\mu \in \mathfrak{g}^\ast$, with respect to which the given
momentum map $\mathbf{J}$ is equivariant. Assume that $G$ acts
freely and properly on $T^\ast Q$, and $\tilde{G}_\mu$ is denoted the
isotropy subgroup of $\mu \in \mathfrak{g}^\ast$ relative to this
affine action $\Theta$, and $\mathcal{O}_\mu= G\cdot \mu
\subset \mathfrak{g}^\ast$ is denoted the G-orbit of the point $\mu$
with respect to the action $\Theta$,
and $\mu$ is a regular value of $\mathbf{J}$.
Then the quotient space $(T^\ast
Q)_\mu=\mathbf{J}^{-1}(\mu)/\tilde{G}_\mu$ is a symplectic
manifold with the symplectic form $\omega_\mu$ uniquely characterized by
$(6.1)$, and the quotient space
$(T^\ast Q)_{\mathcal{O}_\mu}=\mathbf{J}^{-1}(\mathcal{O}_\mu)/ G $
is also a symplectic manifold with the symplectic form
$\omega_{\mathcal{O}_\mu}$ uniquely characterized by $(6.2)$,
see Ortega and Ratiu \cite{orra04}.
Moreover, in this case,
for the given regular point or regular orbit reducible RCH system
$(T^*Q,G,\omega,H,F,W)$, we can also prove the regular point reduction theorem
or regular orbit reduction theorem, by using the above similar ways.
\end{rema}

\section{Some Developments}

In this section, we shall give some generalizations of the above
results from the viewpoint of change of geometrical structures.\\

\noindent{\bf (1) Optimal Reductions of a CH System}\\

It is a natural problem what and how we could do, if we define a controlled
Hamiltonian system on the cotangent bundle $T^*Q$ by using a Poisson
structure, and if symplectic reduction procedure given by Marsden et
al. \cite{mawazh10} does not work or is not efficient enough. In
Wang and Zhang \cite{wazh12}, we study the optimal reduction theory
of a CH system with Poisson structure and symmetry, by using the
optimal momentum map and the reduced Poisson tensor (resp. the reduced
symplectic form). We prove the optimal point reduction,
optimal orbit reduction, and
regular Poisson reduction theorems for the CH system, and explain the
relationships between OpCH-equivalence, OoCH-equivalence,
RPR-CH-equivalence for the optimal reducible CH systems with symmetries
and the CH-equivalence for the associated optimal reduced CH systems.
This paper is published in Jour. Geom. Phys., 62(5)(2012), 953-975.\\

The late Professor Jerrold E. Marsden had joined this research. When the
paper was about to finish, he left us. We are extremely sad. H.Wang
and Z. X. Zhang would like to acknowledge his understanding, support
and help in more than two years of cooperation.\\

\noindent{\bf (2) Singular Reductions of an RCH System}\\

It is worthy of noting that when Lie group $G$ acts only properly on $Q$, does
not act freely, then the reduced space $(T^\ast
Q)_\mu=\mathbf{J}^{-1}(\mu)/G_\mu$ (resp. $(T^\ast
Q)_{\mathcal{O}_\mu}= \mathbf{J}^{-1}(\mathcal{O}_\mu)/G$) is not
necessarily smooth manifold, but just quotient topological space,
and it may be a symplectic Whitney stratified space. In this case,
we study the singular reduction theory of the RCH system
with symmetry, by using the momentum map and the singular reduced
symplectic form in the stratified phase space.
We prove the singular point reduction and
singular orbit reduction theorems for the RCH system, and explain the
relationships between SpCH-equivalence, SoCH-equivalence
for the singular reducible RCH systems with symmetries
and the RCH-equivalence for the associated
singular reduced CH systems. Before Professor Marsden passed away,
we have begun looking for a couple of practical examples which
are RCH-equivalent for singular reduced RCH systems.
But, eight years flied away, we have not found them.
It is not easy to do for the complex structure of stratified phase
space of a singular reduced RCH system.\\

\noindent{\bf (3) Poisson Reduction by
Controllability Distribution for a CH System}\\

It is worthy of noting that when there is no momentum map of Lie group action
for our considered system, then the reduction procedures given in
Marsden et al. \cite{mawazh10} and Wang and Zhang \cite{wazh12} can
not work. One must look for a new way. On the other hand, motivated
by the work of Poisson reductions by distribution for Poisson
manifolds, see Marsden and Ratiu \cite{mara86}, we note that the phase space
$T^*Q$ of the CH system is also a Poisson manifold, and its control
subset $W\subset T^*Q$ is a fiber submanifold. If we assume that
$D\subset TT^*Q |_W$ is a controllability distribution of the CH
system, then we can study naturally the Poisson
reduction by controllability distribution for the CH system. For a
symmetric CH system, and its control subset $W\subset T^*Q$ is a
$G$-invariant fiber submanifold, if we assume that $D\subset TT^*Q
|_W$ is a $G$-invariant controllability distribution of the symmetric CH
system, then we can give Poisson reducible conditions
by controllability distribution for this CH system, and prove the Poisson reducible
property for the CH system and it is kept invariant under the CH-equivalence.
We also study the relationship between Poisson reduction by
$G$-invariant controllability distribution, for the
regular (resp. singular) Poisson reducible CH system, and
Poisson reduction by the reduced controllability distribution, for the associated
reduced CH system. In addition, we can also develop the singular
Poisson reduction and SPR-CH-equivalence for a
CH system with symmetry, and prove the singular Poisson reduction theorem.
See Ratiu and Wang \cite{rawa12} for more details.\\

\noindent{\bf (4) Regular Reduction of a CMH System with
Symmetry of the Heisenberg Group}\\

We consider the regular point reduction of the controlled magnetic Hamiltonian (CMH) system
$(T^\ast Q,\mathcal{H},\omega_Q,H,F,\mathcal{C})$ with symmetry of
the Heisenberg group $\mathcal{H}$. Here the configuration space $Q=
\mathcal{H}\times V, \; \mathcal{H}= \mathbb{R}^2\oplus \mathbb{R},
$ and $V$ is a $k$-dimensional vector space, and the cotangent
bundle $T^*Q$ with the magnetic symplectic form $\omega_Q= \Omega_0-
\pi^*_Q \bar{B},$ where $\Omega_0$ is the usual canonical symplectic
form on $T^*Q$, and $\bar{B}= \pi_1^*B$ is the closed two-form on
$Q$, $B$ is a closed two-form on $\mathcal{H}$ and the projection $\pi_1:
Q=\mathcal{H}\times V \rightarrow \mathcal{H}$ induces the map
$\pi_1^*: T^* \mathcal{H}\rightarrow T^*Q$.\\

A magnetic Hamiltonian system is a Hamiltonian system defined by the
magnetic symplectic form, and a controlled magnetic Hamiltonian (CMH) system on $T^*Q$ is
a magnetic Hamiltonian system $(T^\ast Q,\omega_B,H)$
with external force $F$ and control $W$, where
the magnetic symplectic form $\omega_B= \omega- \pi_Q^*B $ on $T^*Q$,
and $\omega$ is the canonical symplectic form on $T^* Q$
and $B$ is a closed two-form on $Q$, and $F: T^*Q\rightarrow T^*Q$ is
the fiber-preserving map, and $W\subset T^*Q$ is a fiber submanifold.
We note that there is a magnetic term on the cotangent
bundle of the Heisenberg group $\mathcal{H}$, which is related to a
curvature two-form of a mechanical connection determined by the
reduction of center action of the Heisenberg group $\mathcal{H}, $
see Marsden et al. \cite{mamipera98},
such that we can define a controlled magnetic Hamiltonian system with
symmetry of the Heisenberg group $\mathcal{H}$, and study the regular point
reduction of this system. Since a CMH system is also a
RCH system, but its symplectic structure
is given by a magnetic symplectic form. Thus, the set of
the CMH systems is a subset of the set of
the RCH systems, and the set of
the CMH systems with symmetries is also a subset of the set of
the RCH systems with symmetries, and the subset is not complete
under the regular point reduction of the RCH system. It is worthy of noting that
it is different from the regular point reduction of an RCH system
defined on a cotangent bundle with the canonical structure,
the regular point reduction of a CMH system reveals
the deeper relationship of the intrinsic geometrical structures
of the RCH systems on a cotangent bundle. See Wang \cite{wa15a} for more details.

\section{Applications}

In this section, we shall give two applications
for the regular reduction of an RCH system.

\subsection{Two Types of Hamilton-Jacobi Equations }

It is well-known that Hamilton-Jacobi theory is an important research subject
in mathematics and analytical mechanics.
see Abraham and Marsden \cite{abma78}, Arnold
\cite{ar89} and Marsden and Ratiu \cite{mara99},
and the Hamilton-Jacobi equation is
also fundamental in the study of the quantum-classical relationship
in quantization, see Woodhouse \cite{wo92}.
Hamilton-Jacobi theory from the variational
point of view is originally developed by Jacobi in 1866, which states
that the integral of Lagrangian of a system along the solution of
its Euler-Lagrange equation satisfies the Hamilton-Jacobi equation.
The classical description of this problem from the generating function and the geometrical point
of view is given by Abraham and Marsden in \cite{abma78} as follows:
Let $Q$ be a smooth manifold and $TQ$ the tangent bundle, $T^* Q$
the cotangent bundle with the canonical symplectic form $\omega$,
and the projection $\pi_Q: T^* Q \rightarrow Q $ induces the map $
T\pi_{Q}: TT^* Q \rightarrow TQ. $
\begin{theo}
Assume that the triple $(T^*Q,\omega,H)$ is a Hamiltonian system
with Hamiltonian vector field $X_H$, and $W: Q\rightarrow
\mathbb{R}$ is a given generating function. Then the following two assertions
are equivalent:\\
\noindent $(\mathrm{i})$ For every curve $\sigma: \mathbb{R}
\rightarrow Q $ satisfying $\dot{\sigma}(t)= T\pi_Q
(X_H(\mathbf{d}W(\sigma(t))))$, $\forall t\in \mathbb{R}$, then
$\mathbf{d}W \cdot \sigma $ is an integral curve of the Hamiltonian
vector field $X_H$.\\
\noindent $(\mathrm{ii})$ $W$ satisfies the Hamilton-Jacobi equation
$H(q^i,\frac{\partial W}{\partial q^i})=E, $ where $E$ is a
constant.
\end{theo}

From the proof of the above theorem given in
Abraham and Marsden \cite{abma78}, we know that
the assertion $(\mathrm{i})$ with equivalent to
Hamilton-Jacobi equation by the generating function,
gives a geometric constraint condition of the canonical symplectic form
on the cotangent bundle $T^*Q$
for Hamiltonian vector field of the system.
Thus, the Hamilton-Jacobi equation reveals the deeply internal relationships of
the generating function, the canonical symplectic form
and the dynamical vector field of a Hamiltonian system.\\

But, from $\S 2$ and Marsden et al.\cite{mawazh10}, we know that,
the set of Hamiltonian systems with symmetries on a cotangent
bundle is not complete under the Marsden-Weinstein reduction.
Since the symplectic reduced system of a
Hamiltonian system with symmetry defined on the cotangent bundle
$T^*Q$ may not be a Hamiltonian system on a cotangent bundle,
then we cannot give the Hamilton-Jacobi theorem for the Marsden-Weinstein
reduced system just like same as the above Theorem 10.1.
We have to look for a new way. Recently,
Wang in \cite{wa17} prove an important lemma, which is
a modification for the corresponding result of Abraham and Marsden
in \cite{abma78}, such that we can derive precisely
the geometric constraint conditions of
the regular reduced symplectic forms for the
dynamical vector fields of a regular reducible Hamiltonian
system on the cotangent bundle of a configuration manifold,
which are called the Type I and Type II of Hamilton-Jacobi equation,
because they are the development of classical Hamilton-Jacobi equation
given by Abraham and Marsden
in \cite{abma78}, also see Wang \cite{wa17}. \\

Since the Hamilton-Jacobi theory is developed based on the
Hamiltonian picture of dynamics, it is natural idea to extend the
Hamilton-Jacobi theory to the RCH system and its a variety of the reduced systems.
In general, under the actions of external force and control,
an RCH system may not be a Hamiltonian system, and it has yet no generating function,
we cannot give the Hamilton-Jacobi theorems for the RCH system
and its regular reduced systems just like same as the above Theorem 10.1.
But, Wang in \cite{wa13d}derive precisely the geometric constraint conditions of
canonical symplectic form for the
dynamical vector fields of an RCH system on the cotangent
bundle of a configuration manifold, that is,
the Type I and Type II of Hamilton-Jacobi equation,
and  prove that the RCH-equivalence for the RCH systems
leaves the solutions of corresponding Hamilton-Jacobi equations invariant.
Moreover, for the regular point and orbit
reducible RCH systems with symmetries and momentum maps given in $\S 6$,
we also prove two types of Hamilton-Jacobi theorems, and by using the
regular point reduction Theorem 7.2 and the regular orbit reduction Theorem 8.2,
prove the RpCH-equivalence and
RoCH-equivalence for the regular reducible RCH systems leave the
solutions of corresponding Hamilton-Jacobi equations invariant.
In order to describe the impact of different structures of geometry
for the Hamilton-Jacobi theorem, we also prove
two types of Hamilton-Jacobi theorem for a controlled magnetic Hamiltonian (CMH) system,
by using the magnetic symplectic form and the magnetic Hamiltonian vector field.\\

When a Hamiltonian system with nonholonomic constraints,
de Le\'{o}n and Wang in \cite{lewa15} study the
Hamilton-Jacobi theory for
the nonholonomic Hamiltonian system and
the nonholonomic reducible Hamiltonian system
on a cotangent bundle, by using the distributional Hamiltonian system
and the reduced distributional Hamiltonian system. These researches
reveal from the geometrical point of view the internal relationships of
nonholonomic constraints, symplectic forms and
nonholonomic dynamical vector fields of a
mechanical system and its nonholonomic reduced systems.

\subsection{Rigid Body and Heavy Top with Internal Rotors}

Now, it is a natural problem
if there is a practical RCH system and how to show the effect on controls
in regular point reduction and Hamilton-Jacobi theory of the system.
Before, as application examples of theoretical results, Professor
Jerrold E. Marsden gave us two couples of examples in October 2008,
that is, (i) the rigid body with external force torque and that with
internal rotors; (ii) the rigid body with internal rotors and the
heavy top with internal rotors. He told us that these systems are
symmetric, and have $G$-invariant control, may consider reduction
and the RCH-equivalence, but need to calculate them in detail. Thus,
for these examples, we have done three things: (i) to do reduction
by calculation in detail;
(ii) to state the reduced systems to be RCH systems;
(iii) to state every couple of
the reduced systems to be RCH-equivalent.
Because the configuration space of rigid body (resp. heavy
top) is the Lie group $G= \textmd{SO}(3)$ (resp. $G=
\textmd{SE}(3)$), and the configuration space of rigid body with
internal rotors (resp. heavy top with internal rotors) is the
generalized case $G\times V$ of Lie group, where $V$ is a vector space.
In consequence, we first
discuss the regular point reducible RCH systems on a Lie group and
on its generalization, then apply these results to the
concrete examples.
In addition, when we calculated these
examples in detail, we also found the third couple of RCH-equivalent
systems, that is, the rigid body with internal rotors (or external
force torque) and the heavy top, where the heavy top is regarded as
a regular point reducible RCH system without the external force and
control. See Marsden et al. \cite{mawazh10} for more details.
We can also give some other examples, but it is not
easy to give a couple of RCH-equivalent examples, just like same as Professor
Jerrold E. Marsden had done.\\

When we consider the regular point reduction of an RCH system on a Lie
group $G$, the dual of Lie algebra $\mathfrak{g}^\ast$
is a Poisson manifold with respect to
its Lie-Poisson bracket, and the coadjoint orbit $\mathcal{O}_\mu$
through a regular point $\mu\in\mathfrak{g}^\ast$ is the symplectic leaf of
$\mathfrak{g}^\ast. $ We can obtain the $R_p$-reduced RCH system on
$\mathcal{O}_\mu$ by using the $R_p$-reduced symplectic form
$\omega_{\mathcal{O}_\mu}. $ Moreover,
for the regular point reducible RCH system on the generalization of a Lie
group $Q=G\times V$, and $V$ is a vector space,
in this case, for $\mu\in\mathfrak{g}^\ast$, the regular value of
the momentum map $\mathbf{J}_Q $, the $R_p$-reduced phase space is
$\mathcal{O}_\mu \times T^*V\cong \mathcal{O}_\mu \times V\times V^\ast,$ and its
$R_p$-reduced symplectic form is the synthetic of the symplectic form
$\omega_{\mathcal{O}_\mu}$ and the canonical symplectic form
$\omega_{V}$ on $T^*V. $ But, in this case we consider the space
$\mathcal{O}_\mu \times V\times V^\ast$ is the symplectic leaf of
$\mathfrak{g}^\ast \times V\times V^\ast $ with whose Poisson
structure is the synthetic of Lie-Poisson bracket on
$\mathfrak{g}^\ast $ and the Poisson bracket induced from
$\omega_{V}$ on $T^*V. $ In particular,
for the rigid body with three internal rotors, and with a
control torque $u: T^\ast Q \to W $ acting on the rotors, in this case, we take that
$G=\textmd{SO}(3)$ and locally, $V= \mathbb{R}^3, $ then the $R_p$-reduced
controlled rigid body-rotor system is an RCH system on the
symplectic leaf $(\mathcal{O}_\mu \times \mathbb{R}^3 \times
\mathbb{R}^{3*},\tilde{\omega}_{\mathcal{O}_\mu \times \mathbb{R}^3 \times
\mathbb{R}^{3*}}^{-})$. For the heavy top with two pairs of symmetric
internal rotors, and with a control torque $u: T^\ast Q \to W $
acting on the rotors, in this case, we take that
$G=\textmd{SE}(3)= \textmd{SO}(3)\circledS \mathbb{R}^3$
and locally, $V= \mathbb{R}^2, $ by using the Hamiltonian
reduction by stages for semidirect product Lie group,
see Leonard and Marsden \cite{lema97},
for $(\mu,a) \in \mathfrak{se}^\ast(3)$, the regular
value of the momentum map $\mathbf{J}_Q$,
then the $R_p$-reduced controlled heavy top-rotor system
is an RCH system on the symplectic leaf
$(\mathcal{O}_{(\mu,a)} \times \mathbb{R}^2 \times
\mathbb{R}^{2*},\tilde{\omega}_{\mathcal{O}_{(\mu,a)} \times \mathbb{R}^2
\times \mathbb{R}^{2*}}^{-})$.
Thus, we can deal with uniformly the
symplectic reduction of the rigid body, heavy top, as well as them
with internal rotors, such that we can state that all these systems
are the regular point reducible RCH systems and can give their
RCH-equivalences. See Marsden et al. \cite{mawazh10} for more details.\\

On the other hand, we can also get the equations of
rigid body and heavy top as well as them with internal rotors (or
the external force torques), by using Lagrangian and
Euler-Lagrange equation, just like same as that ones have been done in some
references, see Marsden \cite{ma92} and Bloch et al. \cite{blkrmaal92}.
But, it cannot state uniformly that these systems are
the regular point reduced RCH systems, and hence it cannot state
every couple of them to be RCH-equivalent.
Thus, it is very important for a rigorous theoretical work to
offer uniformly a composition of the research results
from a global view point.\\

As an application of the theoretical results for Hamilton-Jacobi theory,
Wang in \cite{wa13d} derive precisely the Type I and Type II of Hamilton-Jacobi equations
for the $R_p$-reduced controlled rigid body-rotor system and
the $R_p$-reduced controlled heavy top-rotor
system on the generalization of rotation group
$\textmd{SO}(3)$ and on the generalization of Euclidean group
$\textmd{SE}(3)$ by calculation in detail, respectively.
In particular, it is worthy of noting that,
the motions of the controlled rigid body-rotor system
and the controlled heavy top-rotor system are different, and
the configuration spaces, the Hamiltonian functions, the actions of Lie groups,
the $R_p$-reduced symplectic forms and the $R_p$-reduced systems of
the controlled rigid body-rotor system
and the controlled heavy top-rotor system, all of them are also different. But,
the two types of Hamilton-Jacobi equations given by calculation in detail
are same, that is, the internal rules are same. This is very important!\\

It is worthy of noting that the method of calculation for the regular point
reductions and the two types of Hamilton-Jacobi equations for the rigid body with
internal rotors and the heavy top with internal rotors is very important and efficient,
and it is generalized and used to the more practical controlled Hamiltonian systems.
For example, Wang in \cite{wa20a, wa13e} applies the research method
to give explicitly the regular point reductions and the two types of Hamilton-Jacobi
equations for the $R_p$-reduced controlled spacecraft-rotor system and
the $R_p$-reduced controlled underwater vehicle-rotor system.\\

Finally, we also note that there have been a
lot of beautiful results of reduction theory of Hamiltonian systems
in celestial mechanics, hydrodynamics and plasma physics. So, it is
an important topic to study the application of symmetric reduction and
Hamilton-Jacobi theory for the controlled Hamiltonian systems
in celestial mechanics, hydrodynamics and
plasma physics. These are our goals in future research.
In particular, it is the key thought of the researches of geometrical mechanics
of the Professor Jerrold E. Marsden to explore and reveal the deeply internal
relationship between the geometrical structure of phase space and the dynamical
vector field of a mechanical system. It is also our goal of pursuing and inheriting.\\

\noindent {\bf Acknowledgments:}
This is a survey to introduce briefly some recent developments
of symmetric reductions of the regular controlled
Hamiltonian systems with symmetries and applications.
The author would like to thank the cultivation and training
of her advisor---Professor Hesheng Hu and to
dedicate the article to her on the occasion of her 90th birthday.\\

Especially grateful to Professor Jerrold E. Marsden,
Professor Tudor S. Ratiu, Professor Manuel de Le\'{o}n, Professor
Arjan Van der Schaft and Professor Juan-Pablo Ortega and MS. Wendy McKay
for their understanding, guiding, support and help
in the study of geometric mechanics.\\

Professor Jerrold E. Marsden is worthy to be
respected. He knew and admitted their wrong and joined us to correct
these wrong. Professor Marsden is a model of excellent
scientists. It is an important task for us to correct and develop well
the research work of Professor Marsden,
such that we never feel sorry for his great cause.\\

Marsden-Weinstein reduction is just like an elegant
cobblestone on the long river bank of history of mankind
civilization, ones share its beauty and happiness.\\


\begin{thebibliography}{30}
\bibitem{abma78}
Abraham R. and Marsden J.E., Foundations of Mechanics, 2$^{nd}$
edition, Addison-Wesley, 1978.
\bibitem{abmara88}
Abraham R., Marsden J.E. and Ratiu T.S., Manifolds, Tensor Analysis
and Applications,  Applied Mathematical Science, \textbf{75},
Springer-Verlag, 1988.
\bibitem{ar89}
Arnold V.I., Mathematical Methods of Classical Mechanics, 2$^{nd}$
edition, Graduate Texts in Mathematics, \textbf{60},
Springer-Verlag, 1989.
\bibitem{blkrmaal92}
Bloch A.M., Krishnaprasad P.S., Marsden J.E., and S\'{a}nchez de
Alvarez G., Stabilization of rigid body dynamics by internal and
external torques, Automatica, \textbf{28}(4), 745--756(1992).
\bibitem{blchlema01}
Bloch A.M., Chang D.E., Leonard N.E. and Marsden J.E., Controlled
Lagrangians and the stabilization of mechanical systems II:
potential shaping,  IEEE Trans, Automatic Control, \textbf{46},
1556--1571(2001).
\bibitem{blle02}
Bloch A.M. and Leonard N.E., Symmetries, conservation laws, and
control, In ``Geometry, Mechanics and Dynamics, Volume in Honour of
the 60th Birthday of J.E. Marsden" (eds. P.Newton, P.Holmes and A.
Weinstein), Springer, New York, 2002.
\bibitem{bllema00}
Bloch A.M., Leonard N.E. and Marsden J.E., Controlled Lagrangian and
the stabilization of mechanical systems I: the first matching
theorem,  IEEE Trans, Automatic Control, \textbf{45},
2253--2270(2000).
\bibitem{bllema01}
Bloch A.M., Leonard N.E. and Marsden J.E., Controlled Lagrangians
and the stabilization of Euler-Poincar$\acute{e}$ mechanical
systems, Int. J. Nonlinear and Robust Control, 11, 191-214(2001).
\bibitem{chbllemawo02}
Chang D.E., Bloch A.M., Leonard N.E., Marsden J.E. and Woolsey C.A.,
The equivalence of controlled Lagrangian and controlled Hamiltonian
systems,  ESAIM Control, Optimisation and Calculus of Variations,
\textbf{8}, 393--422(2002).
\bibitem{chma04}
Chang D.E. and Marsden J.E., Reduction of controlled Lagrangian and
Hamiltonian systems with symmetry, SIAM J. Control Optimization,
\textbf{43}(1), 277--300(2004).
\bibitem{kakost78}
Kazhdan D., Kostant B. and Sternberg S., Hamiltonian group actions dynamical
systems of Calogero type, Comm. Pure Appl. Math. \textbf{31}, 481-508(1978).
\bibitem{lewa15}
de Le\'{o}n M. and Wang H., Hamilton-Jacobi theorems for
nonholonomic reducible Hamiltonian systems on a cotangent bundle,
(arXiv: 1508.07548, a revised version).
\bibitem{lema97}
Leonard N.E, Marsden J.E, Stability and drift of underwater
vehicle dynamics: mechanical systems with rigid motion symmetry,
Physica D, \textbf{105}, 130--162(1997).
\bibitem{lima87}
Libermann P. and Marle C.M., Symplectic Geometry and Analytical
Mechanics, Kluwer Academic Publishers, 1987.
\bibitem{ma76}
Marle C.M., Symplectic manifolds, dynamical groups and Hamiltonian mechanics,
In: Differential Geometry and Relativity, (M. Cahen and M. Flato, eds.),
D. Reidel, Boston, 249-269(1976).
\bibitem{ma92}
Marsden J.E., Lectures on Mechanics, London Mathematical Society
Lecture Notes Series, \textbf{174}, Cambridge University Press,
1992.
\bibitem{mamiorpera07}
Marsden J.E., Misiolek G., Ortega J.P., Perlmutter M. and Ratiu
T.S., Hamiltonian Reduction by Stages, Lecture Notes in Mathematics,
\textbf{1913}, Springer, 2007.
\bibitem{mamipera98}
Marsden J.E., Misiolek G., Perlmutter M. and Ratiu T.S., Symplectic
reduction for semidirect products and central extensions, Diff.
Geom. Appl., \textbf{9}, 173-212 (1998).
\bibitem{mape00}
Marsden J.E. and Perlmutter M., The orbit bundle picture of
cotangent bundle reduction, C. R. Math. Acad. Sci. Soc. R. Can.,
\textbf{22}, 33-54 (2000).
\bibitem{mara86}
Marsden J.E. and Ratiu T.S.: Reduction of Poisson manifolds, Lett.
Math. Phys., \textbf{11(2)}, 161-169, (1986).
\bibitem{mara99}
Marsden J.E. and Ratiu T.S., Introduction to Mechanics and Symmetry,
 2$^{nd}$ edition, Texts in Applied Mathematics, \textbf{17},
Springer-Verlag, 1999.
\bibitem{mamora90}
Marsden J.E., Montgomery R. and Ratiu T.S., Reduction, Symmetry and
Phases in Mechanics, In: Memoirs of the American Mathematical
Society, \textbf{88}, American Mathematical Society, Providence, Rhode
Island, 1990.
\bibitem{mawazh10}
Marsden J.E., Wang H., and Zhang Z.X., Regular reduction of
controlled Hamiltonian system with symplectic structure and
symmetry, Diff. Geom. Appl., \textbf{33}(3), 13-45 (2014), (arXiv:
1202.3564, a revised version).
\bibitem{mawe74}
Marsden J.E. and Weinstein A., Reduction of symplectic manifolds
with symmetry, Rep. Math. Phys., \textbf{5}, 121--130 (1974).
\bibitem{me73}
Meyer K.R., Symmetries and integrals in mechanics, In Peixoto M.
(eds), Dynamical Systems, Academic Press, 259--273 (1973).
\bibitem{nivds90}
Nijmeijer H. and Van der Schaft A.J., Nonlinear Dynamical Control
Systems, Springer-Verlag, 1990.
\bibitem{orra04}
Ortega J.P. and Ratiu T.S., Momentum Maps and Hamiltonian Reduction,
Progress in Mathematics, \textbf{222}, Birkh\"{a}user, 2004.
\bibitem{rawa12}
Ratiu T.S. and Wang H., Poisson reduction by controllability distribution
for a controlled Hamiltonian system, (arXiv: 1312.7047, a revised version).
\bibitem{wa17}
Wang H., Hamilton-Jacobi theorems for regular reducible Hamiltonian
systems on a cotangent bundle, Jour. Geom. Phys., \textbf{119}
82-102, (2017).
\bibitem{wa15a}
Wang H., Regular reduction of a controlled magnetic
Hamiltonian system with symmetry of the Heisenberg group,
(arXiv: 1506.03640, a revised version).
\bibitem{wa13d}
Wang H., Hamilton-Jacobi theorems for regular controlled Hamiltonian
system and its reduced systems, (arXiv: 1305.3457, a revised version).
\bibitem{wa20a}
Wang H., Dynamical equations of the controlled rigid spacecraft with a rotor,
(arXiv: 2005.02221).
\bibitem{wa13e}
Wang H., Symmetric reduction and Hamilton-Jacobi equations of
underwater vehicle-rotor system, (arXiv: 1310.3014, a revised version).
\bibitem{wazh12}
Wang H. and Zhang Z.X., Optimal reduction of controlled Hamiltonian
system with Poisson structure and symmetry, Jour. Geom. Phys.,
\textbf{62}(5), 953-975 (2012).
\bibitem{wo92}
Woodhouse N.M.J. Geometric Quantization, second ed., Clarendon
Press, Oxford, 1992.
\end{thebibliography}
\end{document}